\documentclass[a4,12pt,twoside]{article}
\usepackage{amsfonts}
\usepackage{amssymb}
\usepackage{amsthm}
\usepackage{amsmath}
\usepackage[dvips]{graphicx,color}
\usepackage{latexsym,bm}

\numberwithin{equation}{section}

\def\R2{\mathbb{R}^2}

\begin{document}           

\begin{flushleft}
\renewcommand{\thefootnote}{\fnsymbol{footnote}}
\center{\bf\LARGE  Mean value properties of harmonic functions on
Sierpinski gasket type fractals
}  
\footnotetext{The research of the first author was supported by the
National Science Foundation of China, Grant 10901081, and the
Project Funded by the Priority Academic Program Development of
Jiangsu Higher Education Institutions. }
 \footnotetext{The research of the second author was supported in part by the National Science Foundation, Grant DMS 0652440.}
 \vskip0.5cm \center{\large
 HUA QIU} AND \large{ROBERT S. STRICHARTZ\hspace{0.1cm} }

\end{flushleft}

\noindent{\bf Abstract.} In this paper, we establish an analogue of
the classical mean value property for both the  harmonic functions
and some general functions in the domain of the Laplacian on the
Sierpinski gasket. Furthermore, we extend the result to some other
p.c.f. fractals with Dihedral-$3$ symmetry.

\noindent\textbf{Keywords.} Sierpinski gasket, Laplacian, harmonic
function, mean value property, analysis on fractals.

\noindent\textbf{Mathematics Subject Classification (2000)} 28A80

\hspace{0.2cm}
\\

\section{Introduction}

It is well known that harmonic functions ($i.e.$, solutions of the
Laplace equation $\Delta u=0$, where
$\Delta=\sum_{i=1}^d\frac{\partial^2}{\partial x_i^2}$) possess the
\emph{mean value property}: Namely, if $u$ is harmonic on a domain
$\Omega\subset \mathbb{R}^d$, then for every closed ball
$B_r(x)\subset\Omega$  of a center $x\in \Omega$ and radius $r>0$
the average of $u$ over $B_r(x)$ equals to the value of $x$, i.e.,
$$\frac{1}{|B_r(x)|}\int_{B_r(x)}u(y)dy=u(x),$$ where $|B_r(x)|$
is the volume of the ball $B_r(x)$. There is a similar statement for
mean values on spheres. More generally, if $u$ is not assumed
harmonic but $\Delta u$ is a continuous function, then
\begin{equation}\label{0}
\lim_{r\rightarrow
0}\frac{1}{r^2}\left(\frac{1}{|B_r(x)|}\int_{B_r(x)}u(y)dy-u(x)\right)=c_n\Delta
u(x)
\end{equation} for the appropriate dimensional constant $c_n$.

What are the fractal analogs of these results? The analytic theory
on p.c.f. fractals was developed by Kigami $\cite{Ki1,Ki2,Ki3}$
following the work of several probabilists who constructed
stochastic processes analogous to Brownian motion, thus obtaining a
Laplacian indirectly as the generator of the process. See the book
of Barlow $\cite{Ba}$ for an account of this development. Since
analysis on fractals has been made possible by the analytic
definition of Laplacian, it is natural to explore the properties of
these fractal Laplacians that are natural analogs of results that
are known for the usual Laplacian. As for the fractal analog of the
mean value property, we won't state the nature of the sets on which
we do the averaging here, but will say that if $K$ is a fractal set
and $x\in K$, we investigate whether there is a sequence of sets
$B_k(x)$ containing $x$ with $\bigcap_kB_k(x)=\{x\}$ such that
$$\frac{1}{\mu(B_k(x))}\int_{B_k(x)}u(y)dy=u(x)$$ for every harmonic
function $u$. Moreover, for general $u$ not assumed harmonic, is
there a formula analogous to $(\ref{0})$?

In the present paper, we will mainly deal with the Sierpinski gasket
$\mathcal{SG}$. This set is a key example of fractals on which a
well established theory of Laplacian exists $[3-7]$. Since the mean
value property plays a very important role in the usual theory of
harmonic functions, it is of independent interest to understand the
similar property of harmonic functions on the Sierpinski gasket. We
will prove that for each point $x\in\mathcal{SG}\setminus V_0$,
($V_0$ is the boundary of $\mathcal{SG}$.) there is a sequence of
\emph{mean value neighborhoods} $B_k(x)$ depending only on the
location of $x$ in $\mathcal{SG}$. $\{B_k(x)\}$ forms a system of
neighborhoods of the point $x$ satisfying $\bigcap_{k}B_k(x)=\{x\}$.
On such sequences, we get the fractal analogs of the  mean value
properties of both the harmonic functions and the general functions
which belong to the domain of the fractal Laplacian satisfying some
natural continuity assumption. We also investigate the extent to
which our method can be applicable to other p.c.f. self-similar
sets, but it seems that it strongly depends on the symmetric
properties of both the geometric structure and the harmonic
structure of the fractals.

The paper is organized as follows: In Section 2 we briefly introduce
some key notions from analysis on the Sierpinski gasket. In Section
3 and Section 4, we prove the mean value property for harmonic
functions and general functions on $\mathcal{SG}$ respectively.
Section 5 contains a further  extension of the mean value property
to p.c.f. self-similar fractals with Dihedral-3 symmetry. An
interesting open question is to what extent the results of Section 4
can be extended to this class of fractals. See $\cite{Blank}$ for a
related result concerning solutions of divergence form elliptic
operators.

\section{Analysis on the Sierpinski gasket}

For the convenience of the reader, we collect some key facts from
analysis on $\mathcal{SG}$ that we need to state and prove our
results. These come from Kigami's theory of analysis on fractals,
and may be found in $\cite{Ki1,Ki2,Ki3}$. An elementary exposition
may be found in $\cite{Str1,Str2}$.
  Recall that $\mathcal{SG}$ is the
attractor of the i.f.s (iterated function system) in the plane
consisting of three homotheties $\{F_0, F_1, F_2\}$ with contraction
ratio $1/2$ and fixed points equal to the three vertices $\{q_0,
q_1, q_2\}$ of an equilateral triangle. Then $\mathcal{SG}$ is the
unique nonempty compact set satisfying
\begin{equation}\label{1}
\mathcal{SG}=\bigcup_{i=0}^2F_i(\mathcal{SG}).
\end{equation}
 We refer to the sets $F_i(\mathcal{SG})$ as \emph{cells} of level one, and by
iterating $(\ref{1})$ we obtain the splitting of $\mathcal{SG}$ into
cells of higher level. For a word $w=(w_1,w_2,\cdots,w_m)$ of length
$m$, the set $F_w(\mathcal{SG})=F_{w_1}\circ F_{w_2}\circ\cdots\circ
F_{w_m}(\mathcal{SG})$ with $w_i\in\{0,1,2\}$, is called an
$m$-cell. The fractal $\mathcal{SG}$ can be realized as the limit of
a sequence of graphs $\Gamma_0, \Gamma_1,\cdots$ with vertices
$V_0\subseteq V_1\subseteq\cdots$. The initial graph $\Gamma_0$ is
just the complete graph on $V_0=\{q_0,q_1,q_2\}$,  which is
considered the boundary of $\mathcal{SG}$. See Fig. 2.1. Note that
$\mathcal{SG}$ is connected, but just barely: there is a dense set
of points $\mathcal{J}$, called \emph{junction points}, defined by
the condition that $x\in \mathcal{J}$ if and only if
$U\setminus\{x\}$ is disconnected for all sufficiently small
neighborhoods $U$ of $x$. It is easy to see that $\mathcal{J}$
consists of all images of $\{q_0, q_1, q_2\}$ under iterates of the
i.f.s. The vertices $\{q_0, q_1, q_2\}$ are not junction points. All
other points in $\mathcal{SG}$ will be called \emph{generic points}.
In the $\mathcal{SG}$ case, $\mathcal{J}=V_*\setminus V_0$, where
$V_*=\bigcup_m V_m$. However, it is not true for general p.c.f.
self-similar sets.
\begin{figure}[h]
\begin{center}
\includegraphics[width=13cm,totalheight=3.6cm]{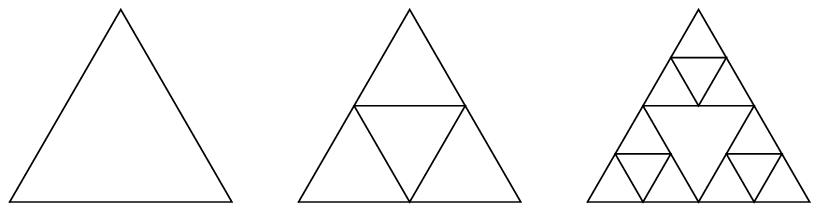}
\setlength{\unitlength}{1cm}
\begin{picture}(0,0) \thicklines
\put(-11.4,-0.2){$\Gamma_{0}$}
\put(-6.9,-0.2){$\Gamma_{1}$}\put(-2.4,-0.2){$\Gamma_2$}
\end{picture}
\begin{center}
\textbf{Fig. 2.1.} \small{The first $3$ graphs, $\Gamma_0, \Gamma_1,
\Gamma_2$ in the approximation to the Sierpinski gasket.}
\end{center}
\end{center}
\end{figure}
In all that follows, we assume that $\mathcal{SG}$ is equipped with
the self-similar probability measure $\mu$ that assigns the measure
$3^{-m}$ to each $m$-cell.

We define the unrenormalized energy of a function $u$ on $\Gamma_m$
by
$$E_m(u)=\sum_{x\sim_m y}(u(x)-u(y))^2.$$
The energy renormalization factor is $r=\frac{3}{5}$, so the
renormalized graph energy on $\Gamma_m$ is
$$\mathcal{E}_m(u)=r^{-m}E_m(u),$$ and we can define the \emph{fractal
energy} $\mathcal{E}(u)=\lim_{m\rightarrow\infty}\mathcal{E}_m(u)$.
We define $dom\mathcal{E}$ as the space of continuous functions with
finite energy. Then $\mathcal{E}$ extends by polarization to a
bilinear form $\mathcal{E}(u,v)$ which serves as an inner product in
this space.

The standard Laplacian may then be defined using the weak
formulation: $u\in dom\Delta$ with $\Delta u=f$ if $f$ is
continuous, $u\in dom\mathcal{E}$, and
\begin{equation*}
\mathcal{E}(u,v)=-\int fvd\mu
\end{equation*}
 for all $v\in
dom_0\mathcal{E}$, where $dom_0\mathcal{E}=\{v\in \mathcal{E}:
v|_{V_0}=0\}$. There is also a pointwise formula (which is proven to
be equivalent in $\cite{Str2}$) which, for points in $V_*\setminus
V_0$ computes
$$\Delta u(x)=\frac{3}{2}\lim_{m\rightarrow\infty}5^m\Delta_m u(x),$$
where $\Delta_m$ is a discrete Laplacian associated to the graph
$\Gamma_m$, defined by
$$\Delta_m u(x)=\sum_{y\sim_m x}(u(y)-u(x))$$ for $x$ not on
the boundary.

It is not necessary to invoke the measure to define \emph{harmonic
functions}, although it is true that these are just the solutions of
$\Delta h=0$. The more direct definition is that
$$h(x)=\frac{1}{4}\sum_{y\sim_m x}h(y)$$ for every nonboundary
point and every $m$. This can be viewed as a mean value property of
$h$ at the junction points. The space of harmonic functions is
$3$-dimensional and the values at the $3$ boundary points may be
freely assigned. Moreover, there is a simple efficient algorithm,
the $``\frac{1}{5}-\frac{2}{5}$ rule", for computing the values of a
harmonic function exactly at all vertex points in terms of the
boundary values. The harmonic functions satisfy the \emph{maximum
principle}, i.e., the maximum and minimum are attained on the
boundary and only on the boundary if the function is not constant.
We call a continuous function $h$ a \emph{piecewise harmonic spline}
of level $m$ if $h\circ F_w$ is harmonic for all $|w|=m$.

The Laplacian satisfies the scaling property
$$\Delta(u\circ F_i)=\frac{1}{5}(\Delta u)\circ F_i$$
and by iteration
$$\Delta(u\circ F_w)=\frac{1}{5^m}(\Delta u)\circ F_w$$ for $F_w=F_{w_1}\circ F_{w_2}\circ\cdots\circ
F_{w_m}$.

Although there is no satisfactory analogue of gradient, there is a
\emph{normal derivative} $\partial_n u(q_i)$ defined at boundary
points by
$$\partial_n u(q_i)=\lim_{m\rightarrow\infty}\sum_{y\sim_m q_i}r^{-m}(u(q_i)-u(y)),$$
the limit existing for all $u\in dom\Delta$. The definition may be
localized to boundary points of cells: for each point $x\in
V_m\setminus V_0$, there are two cells containing $x$ as a boundary
point, hence two normal derivatives at $x$. For $u\in dom\Delta$,
the normal derivatives at $x$ satisfy the \emph{matching condition}
that their sum is zero. The matching conditions allow us to glue
together local solutions to $\Delta u=f$.

As is shown in $\cite{Ki1,Ki2,Str2}$, the Dirichlet problem for the
Laplacian can be solved by integrating against an explicitly given
\emph{Green's function}. Recall that the Green's function $G(x,y)$
is a uniform limit of $G_M(x,y)$ as $M$ goes to the infinity, with
$G_M$ defined by
$$G_M(x,y)=\sum_{m=0}^M\sum_{z,z'\in V_{m+1}\setminus
V_m}g(z,z')\psi_{z}^{(m+1)}(x)\psi_{z'}^{(m+1)}(y)$$ and
\begin{eqnarray*}
\left\{\begin{array}{l}
  g(z,z)=\frac{9}{50} r^m \mbox{ for } z\in V_{m+1}\setminus V_m, \\
  g(z,z')=\frac{3}{50}r^m \mbox{ for } z,z'\in V_{m+1}\setminus V_m
  \mbox{ with } z,z'\in F_w(\mathcal{SG}) \mbox{ for } |w|=m, \mbox{ and }
  z\neq z',
\end{array}\right.
\end{eqnarray*}
 where $\psi_z^{m}(x)$ denotes a piecewise harmonic spline of level $m$ satisfying $\psi_{z}^{(m)}(x)=\delta_{z}(x)$ for $x\in V_m$.

\section{Mean value property of harmonic functions on $\mathcal{SG}$}

\textbf{Lemma 3.1.} \emph{(a) Let $C$ be any cell with boundary points
$p_0,p_1,p_2$, and $h$ any harmonic function. Then
$$\frac{1}{\mu(C)}\int_C hd\mu=\frac{1}{3}(h(p_0)+h(p_1)+h(p_2)).$$}

\emph{(b) Let $p$ be any junction point, and $C_1$, $C_2$ the two
$m$-cells containing $p$. Then
$$\frac{1}{\mu(C_1\cup C_2)}\int_{C_1\cup C_2}hd\mu=h(p).$$}

\emph{Proof.} The space of harmonic functions on $C$ is
three-dimensional. A simple basis $\{h_0,h_1,h_2\}$ is obtained by
taking $h_j(p_j)=1$ and $h_j(p_k)=0$ for $k\neq j$. Noticing that
$h_0+h_1+h_2$ is identically $1$ on $C$, by symmetry, $\int_C h_i
d\mu=\frac{1}{3}\mu(C)$ for each $i$. Hence (a) follows. (b) follows
by combining (a) for $C=C_1$ and $C=C_2$ with the mean value
property of $h$ at $p$. $\Box$

Note that (b) gives a trivial solution to the problem of finding
mean value neighborhoods for junction points.

\begin{figure}[h]
\begin{center}
\includegraphics[width=9.5cm,totalheight=4.5cm]{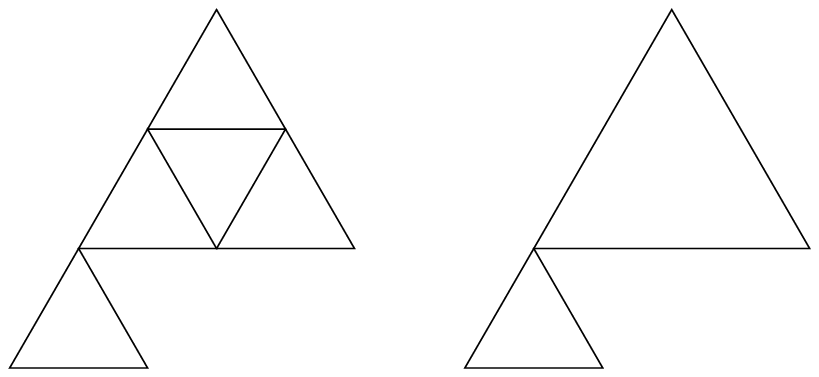}
\setlength{\unitlength}{1cm}
\begin{picture}(0,0) \thicklines
\put(-9,0.5){$C_1$} \put(-7,0.4){$D_w$}\put(-4,0.4){$D_w^{(1)}$}
\put(-1.8,0.4){$D_w$} \put(-9.25,1.5){$p_1$} \put(-7.3,1.3){$p_2$}
\put(-8.2,1.9){$C_w$} \put(-6.6,1.9){$C_2$} \put(-7.4,3.3){$C_0$}
\put(-8.42,2.9){$p_0$} \put(-2.3,2.4){$D_w^{(2)}$}
\end{picture}
\begin{center}
\textbf{Fig. 3.1.} \small{$C_w$ and its three neighboring cells.

The right part of the figure refers to the proof of Lemma 4.1.}
\end{center}
\end{center}
\end{figure}
Given a point $x$ in $\mathcal{SG}\setminus V_0$, consider any cell
$F_w(\mathcal{SG})$(denote it by $C_w$) containing the point $x$,
with boundary points $F_{w}q_i=p_i$. Choose the cell $C_w$ small
enough, such that it does not intersect $V_0$. Then it must have
three neighboring cells $C_0$, $C_1$ and $C_2$ of the same level
with $C_i$ intersecting $C_w$ at $p_i$. Denote by $D_w$ the union of
$C_w$ and its three neighbors. See Fig. 3.1. In this section, we
will describe a method to find a subset $B$ of $D_w$, containing
$C_w$, such that for any harmonic function $h$, the mean value of
$h$ over $B$ is equal to its value at $x$, i.e., $M_B(h)=h(x)$ where
$M_B(h)$ is defined by
$$M_B(h)=\frac{1}{\mu(B)}\int_B hd\mu.$$ Then we will call the set $B$  a \textit{$k$
level mean value neighborhood} of $x$ associated to  $C_w$ where $k$
is the length of $w$.

Let $h$ be a harmonic function on $\mathcal{SG}$. The harmonic
extension algorithm implies that there exist coefficients
$\{a_i(x)\}$ depending only on the relative position of $x$ and
$C_w$ such that
$$h(x)=\sum_i a_i(x)h(p_i).$$ Moreover, since constants are harmonic
we must have $$\sum_i a_i(x)=1$$ and by the maximum principle all
$a_i(x)\geq 0$. Let $W$ denote the triangle in $\mathbb{R}^3$ with
boundary points $(1,0,0), (0,1,0)$ and $(0,0,1)$ and $\pi_W$ the
plane in $\mathbb{R}^3$ containing $W$.  So
$\{(a_0(x),a_1(x),a_2(x))\}\in W$ for any $x\in C_w$. However, not
every point in $W$ occurs in this way.

 On the other hand,  given a set $B$
such that $C_w\subset B\subset D_w$, by linearity we have
\begin{equation}\label{5}
M_B(h)=\sum_i a_i h(p_i)
\end{equation}
 for some coefficients
$(a_0,a_1,a_2)$ depending only on the relative geometry of $B$ and
$C_w$. Again we must have $\sum a_i=1$ by considering $h\equiv 1$.
So $(a_0,a_1,a_2)\in\pi_W$. (Later we will show that $(a_0,a_1,a_2)$
does not have to belong to $W$ for some sets $B$.) Thus we have a
map, denoted by $\mathcal{T}$ from the collection of $B$'s to
$\pi_W$. If we can show that the image of the map $\mathcal{T}$
covers the triangle $W$ for some reasonable class of sets $B$, then
we can get a set $B$ over which
 the mean value property holds for all harmonic functions. Moreover,
 if we can prove $\mathcal{T}$ is one-to-one, then we get a
  mean value neighborhood $B$ of $x$ associated to $C_w$, that is
  unique within the collection of sets we are considering.

The above is the basic idea of our method. Hence, the remaining task
in this section is to find a suitable class $\mathcal{B}$ of sets
$B$  such that there is a map $\mathcal{T}$ from $\mathcal{B}$ to
$\pi_W$, such that $\mathcal{T}(\mathcal{B})$ covers the triangle
$W$. Comparing with the usual mean value neighborhoods (they are
just balls in the Euclidean case), it is reasonable to require $B$
to be as simple as possible. They should be connected, possess some
symmetry properties, depend only on the relative geometry of $x$ and
$C_w$, and be independent of the level of $C_w$ and the location of
$C_w$.

In the following, we use $\rho$ to denote the distance from $p_0$ to
the line containing $p_1$ and $p_2$, namely, $\rho$ is the length of
the height of the minimal equilateral triangle containing $C_w$.
Call $\rho$ the \emph{size} of $C_w$.

 \textbf{Definition 3.1.} \emph{Let
$c_0,c_1,c_2$ be three real numbers satisfying $0\leq c_i\leq 1$,
denote by $B(c_0,c_1,c_2)$ the set
$$B(c_0,c_1,c_2)=C_w\cup
E_0\cup E_1\cup E_2,$$ where each $E_i$ is a sub-triangle domain in
$C_i$ obtained by cutting $C_i$ symmetrically with a line at the
distance $c_i\rho$ away from the vertex $p_i$.}

\textbf{Remark.} \emph{See Fig. 3.2 for a sketch of
$B(c_0,c_1,c_2)$. For example, $B(0,0,0)=C_w$ and $B(1,1,1)=D_w$.
Denote by
$$\mathcal{B}=\{B(c_0,c_1,c_2): 0\leq c_i\leq 1\}$$ the natural $3$-parameter family of
all such sets. Each member of $\mathcal{B}$ contains $C_w$ and is
contained in $D_w$. Denote by $$\sigma: \mathcal{B} \mapsto \Lambda
$$ the natural one-to-one projection with
$\sigma(B(c_0,c_1,c_2))=(c_0,c_1,c_2)$, where
$\Lambda=\{(c_0,c_1,c_2): 0\leq c_i\leq 1\}$.}

\begin{figure}[h]
\begin{center}
\includegraphics[width=6.2cm,totalheight=6.4cm]{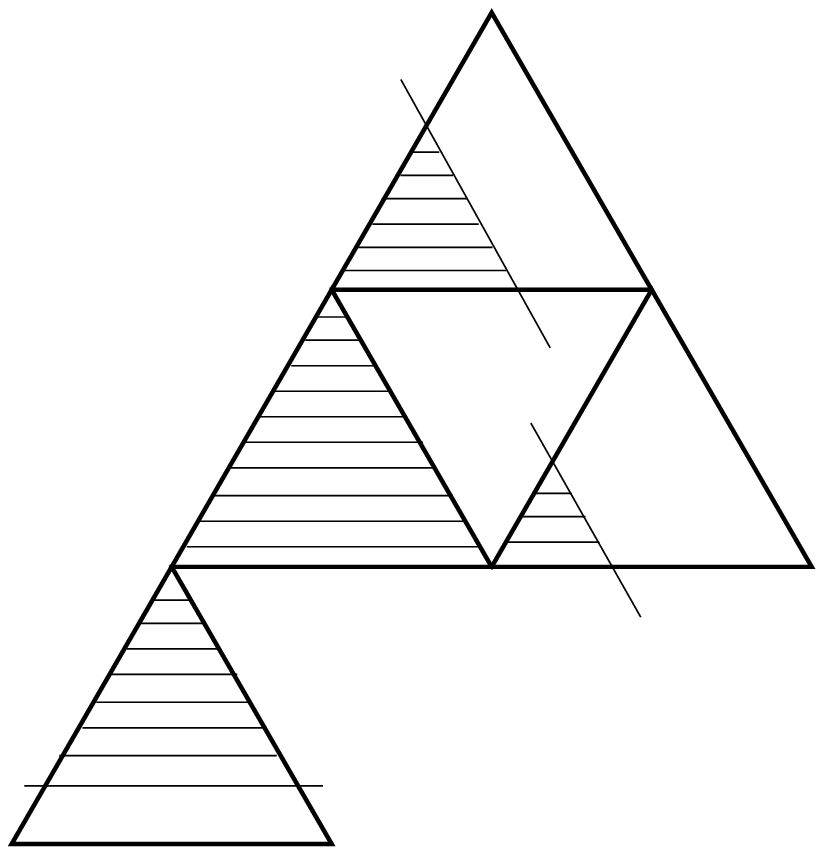}
\setlength{\unitlength}{1cm}
\begin{picture}(0,0) \thicklines
\put(-4.3,2.9){$C_w$} \put(-3.5,4.5){$E_0$} \put(-2.55,2.3){$E_2$}
\put(-4.4,4.3){$p_0$}\put(-5.5,2.3){$p_1$}\put(-3.0,1.9){$p_2$}\put(-5.4,1.0){$E_1$}\put(-2.4,1.0){$B(c_0,c_1,c_2)$}
\end{picture}
\begin{center}
\textbf{Fig. 3.2.} \small{The relative geometry of $B(c_0,c_1,c_2)$
and $C_w$.}
\end{center}
\end{center}
\end{figure}

For each vector $(c_0,c_1,c_2)\in \Lambda$, there is a unique vector
$(a_0,a_1,a_2)\in \pi_W$ corresponding to the set $B(c_0,c_1,c_2)$,
satisfying $(\ref{5})$ where $B$ is replaced by $B(c_0,c_1,c_2)$.
This defines a map $T$ from $\Lambda$ to $\pi_W$. Then $\mathcal{T}$
described above from $\mathcal{B}$ to $\pi_W$ is exactly
$\mathcal{T}=T\circ \sigma$.

 The following lemma shows that the value
$T(c_0,c_1,c_2)$ is independent of the particular choice of $C_w$,
which benefits from the symmetric properties of the set
$B(c_0,c_1,c_2)$.

\textbf{Lemma 3.2.}\emph{ $T(c_0,c_1,c_2)$ is independent of the
particular choice of $C_w$.}

\textit{Proof.} Let $h$ be a harmonic function. First we consider
the integral $\int_{E_i} hd\mu$. Denote by $\{s_i,t_i,p_i\}$ the
boundary points of $C_i$. By linearity,
$\frac{1}{\mu(C_w)}\int_{E_i} hd\mu$ can be expressed as a
non-negative linear combination of $\{h(s_i),h(t_i),h(p_i)\}$, which
by symmetry must have the form
\begin{equation}\label{2}
\int_{E_i} hd\mu=\left(m_ih(p_i)+n_i(h(s_i)+h(t_i))\right)\mu(C_w),
\end{equation} for some appropriate non-negative
coefficients $m_i,n_i$. Notice that in $(\ref{2})$, the coefficients
$m_i,n_i$ are independent of the location of $C_i$ in
$\mathcal{SG}$. Actually, they only depend on the relative position
of $E_i$ in $C_i$, i.e., $m_i,n_i$ depend only on $c_i$. Using the
mean value property at $p_i$, namely
$$4h(p_i)=h(p_{i-1})+h(p_{i+1})+h(s_i)+h(t_i),$$ we obtain
\begin{eqnarray*}
\int_{E_i}hd\mu&=&(m_ih(p_i)+n_i(4h(p_i)-h(p_{i-1})-h(p_{i+1})))\mu(C_w)\\&=&((4n_i+m_i)h(p_i)-n_i(h(p_{i-1})+h(p_{i+1})))\mu(C_w).
\end{eqnarray*}

Notice that the ratio of $\mu(E_i)$ to $\mu(C_i)$ also depends only
on $c_i$. Combined with Lemma 3.1(a),  we see that
$(a_0,a_1,a_2)=T(c_0,c_1,c_2)$ is independent of the particular
choice of $C_w$, depending only on $(c_0,c_1,c_2)$.$\Box$

We will show the image of the map $\mathcal{T}$ covers the triangle
$W$. More precisely, $T(c_0,c_1,c_2)$ will fill out a set
$\widetilde{W}$ which is a bit larger than $W$. Denote by
$P_0=(1,0,0)$, $P_1=(0,1,0)$ and $P_2=(0,0,1)$ the three boundary
points of the triangle $W$ in $\mathbb{R}^3$ and by $O$ the center
point of $W$.
\begin{figure}[h]
\begin{center}
\includegraphics[width=4.2cm,totalheight=4.1cm]{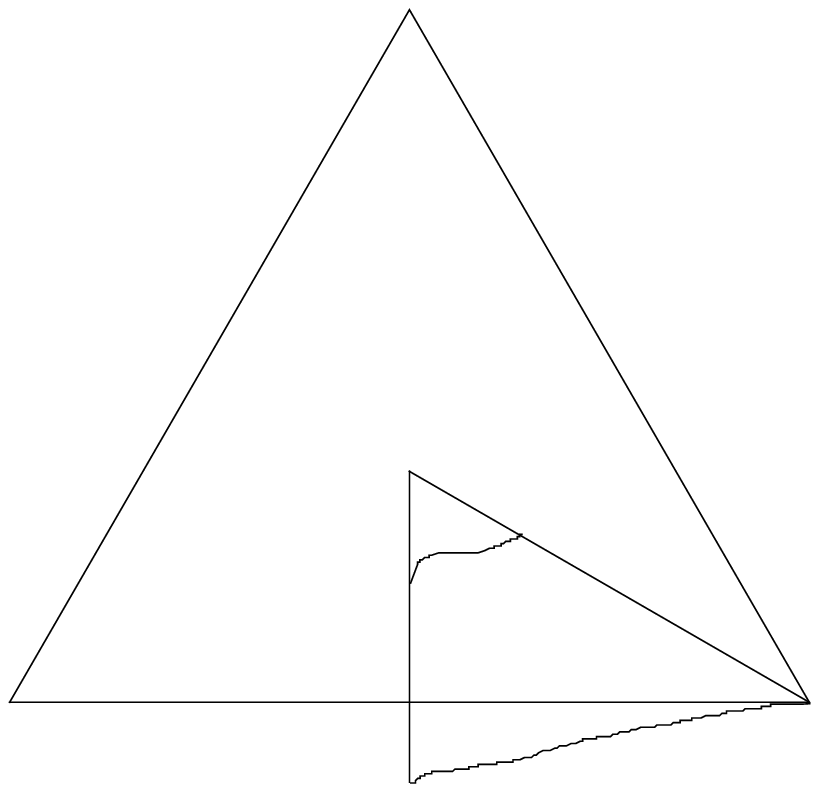}
\setlength{\unitlength}{1cm}
\begin{picture}(0,0) \thicklines
\put(-2.6,1.7){$O$}\put(-2.6,-.3){$Q_0$}\put(-2.1,0.8){$\Gamma_y$}\put(-1.3,0){$\Gamma_1$}\put(-0.3,0.3){$P_2$}
\end{picture}
\begin{center}
\textbf{Fig. 3.3.} \small{a $1/6$ region of $\widetilde{W}$
surrounded by $\overline{OQ_0}$, $\overline{OP_2}$ and
$\widehat{P_2Q_0}$.}
\end{center}
\end{center}
\end{figure}

\textbf{Lemma 3.3.} \emph{$T(0,0,1)=P_2$ and $T(0,1,1)=Q_0$ where
$Q_0=\{-\frac{1}{9},\frac{5}{9},\frac{5}{9}\}$ is a point in $\pi_W$
located outside of $W$.}

\textit{Proof.} From Definition 3.1, $B(0,0,1)=C_w\cup C_2$. Hence
by Lemma 3.1(b), for any harmonic function $h$, we have
$M_{B(0,0,1)}(h)=h(p_2)$. This implies $T(0,0,1)=P_2$. Similarly,
$B(0,1,1)=C_w\cup C_1\cup C_2$, then for any harmonic function $h$,
still using Lemma 3.1, we get
\begin{eqnarray*}
M_{B(0,1,1)}(h)&=&\frac{1}{3\mu(C)}\left(\int_{C_w\cup
C_1}hd\mu+\int_{C_w\cup
C_2}hd\mu-\int_{C_w}hd\mu\right)\\&=&-\frac{1}{9}h(p_0)+\frac{5}{9}h(p_1)+\frac{5}{9}h(p_2),
\end{eqnarray*}
which gives $T(0,1,1)=Q_0$. $\Box$

\textbf{Lemma 3.4.}\emph{ $T(\{(0,c,1): 0\leq c\leq 1\})$ is a
continuous curve lying outside of $W$, joining $P_2$ and $Q_0$. (See
Fig. 3.3.)}

\emph{Proof.} From Lemma 3.3, by varying $c$ continuously between
$0$ and $1$ we trace a continuous curve $\widehat{P_2Q_0}$ joining
$P_2$ and $Q_0$. So we only need to prove the curve
$\widehat{P_2Q_0}$ lies outside of $W$. To prove this, we consider
the set $B=B(0,c,1)$ for $0\leq c\leq 1$. In this case $$B=C_w\cup
E_1\cup C_2.$$

Given a harmonic function $h$, by the proof of Lemma 3.2, we have
\begin{equation*}
\int_{E_1} hd\mu=((4n_1+m_1)h(p_1)-n_1(h(p_{0})+h(p_{2})))\mu(C_w),
\end{equation*} for some appropriate non-negative
coefficients  $m_1,n_1$ depending only on $c$.

On the other hand, we have
$$\int_{C_w\cup C_2}hd\mu=2h(p_2)\mu(C_w),$$
by Lemma 3.1(b).

Hence \begin{eqnarray*}\int_Bhd\mu&=&\int_{E_1}hd\mu+\int_{C_w\cup
C_2}hd\mu\\
&=&(-n_1h(p_0)+(4n_1+m_1)h(p_1)+(2-n_1)h(p_2))\mu(C_w).
\end{eqnarray*}

The coefficient of $h(p_0)$ is always less than $0$. Moreover, it
equals to $0$ if and only if $E_1=\emptyset$ (c=0). Hence $T(0,c,1)$
will always lie on the outside of the triangle $W$ as $c$ varies
between $0$ and $1$. $\Box$

Now we come to the main result of this section.

\textbf{Theorem 3.1.} \emph{The map $\mathcal{T}$ from $\mathcal{B}$
to $\pi_W$ fills out a region $\widetilde{W}$ which contains the
triangle $W$.}

\textit{Proof.} We only need to prove the map $\mathcal{T}$ from
$\mathcal{B}$ to $\pi_W$ fills out a 1/6 region surrounded by  the
line segments $\overline{OQ_0}$, $\overline{OP_2}$ and the curve
$\widehat{P_2Q_0}$ as shown in Fig. 3.3. Then we will get the
desired result by exploiting the symmetry.

Consider a subfamily $\mathcal{B}_1=\{B(0,0,c):0\leq c\leq 1\}$ of
$\mathcal{B}$. If we restrict the map $\mathcal{T}$ to
$\mathcal{B}_1$, by varying $c$ continuously between $0$ and $1$ we
trace a curve (it is a line segment, which follows from the symmetry
of $E_2$) in $W$ joining the center $O$ and the vertex point $P_2$.

Consider another subfamily $\mathcal{B}_2=\{B(0,c,c):0\leq c\leq 1
\}$ of $\mathcal{B}$. If we restrict the map $\mathcal{T}$ to
$\mathcal{B}_2$, by varying $c$ continuously between $0$ and $1$ we
trace a curve (it is also a line segment, which follows from the
symmetric effect of $E_1$ and $E_2$) in $W$ joining the center $O$
and the  point $Q_0$ across the boundary line $\overline{P_1P_2}$
with $Q_0$ located outside of $W$, where $Q_0$ is the point defined
in Lemma 3.3.

Fix a number $0\leq y\leq 1$. Consider a subfamily
$\mathcal{C}_y=\{B(0,c,y):0\leq c\leq y\}$ of $\mathcal{B}$. If we
restrict the map $\mathcal{T}$ to $\mathcal{C}_y$, by varying $c$
continuously between $0$ and $y$ we trace a curve $\Gamma_y$ joining
the two points $T(0,0,y)$ and $T(0,y,y)$. The first endpoint
$T(0,0,y)$ lies on the line segment $\overline{OP_2}$ and the second
endpoint $T(0,y,y)$ lies on the line segment $\overline{OQ_0}$. (See
Fig. 3.3. for $\Gamma_y$.) When $y=0$, the curve $\Gamma_0$ draws
back to the single center point $O$. When $y=1,$ by Lemma 3.4, the
curve $\Gamma_1$ is a continuous curve located outside of the
triangle $W$.  Moreover, $P_2$ is the only common points of
$\Gamma_1$ and $W$. Hence if we vary $y$ continuously  between $0$
and $1$, we can fill out the $1/6$ region surrounded by  the line
segments $\overline{OQ_0}$, $\overline{OP_2}$ and the curve
$\widehat{P_2Q_0}$. $\Box$

\textbf{Remark.}\emph{ In the proof of the above theorem, we
actually only consider those sets $B$ in $\mathcal{B}$ which are
contained in the union of $C_w$ and subsets of only two neighbors.
See Fig. 3.4. Of course, the map $\mathcal{T}$ restricted to this
subfamily is one-to-one, which can be easily seen from the proof.
Hence instead of $\mathcal{B}$, the map $\mathcal{T}$ is one-to-one
from $\mathcal{B}^*$ onto $\widetilde{W}$, where
$$\mathcal{B}^*=\{B(0,c_1,c_2):0\leq c_i\leq 1\}\cup\{B(c_0,0,c_2):0\leq c_i \leq 1\}\cup \{B(c_0,c_1,0): 0\leq c_i\leq
1\}.$$}
\begin{figure}[h]
\begin{center}
\includegraphics[width=11cm,totalheight=4.1cm]{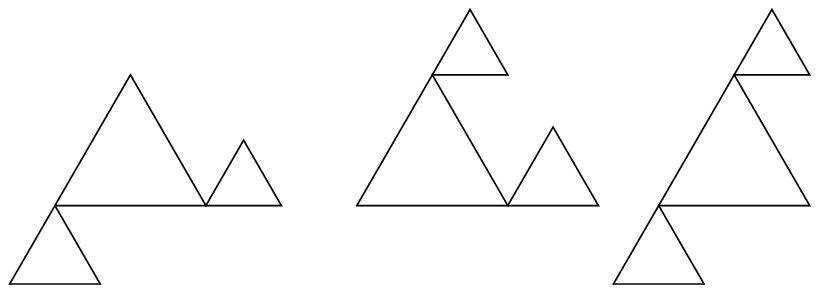}
\begin{center}
\textbf{Fig. 3.4.} \small{The 3 shapes of $B\in \mathcal{B}^*$
associated to $C_w$ shown in Fig. 3.1.}
\end{center}
\end{center}
\end{figure}

Based on the discussion in the beginning of this section, we then
have

\textbf{Theorem 3.2.}\emph{ For each point $x\in
\mathcal{SG}\setminus V_0$, there exists a system of mean value
neighborhoods $B_k(x)$ with $\bigcap_k B_k(x)=\{x\}$.}

\textit{Proof.} Let $k_0$ be the smallest value of $k$ such that
there exists a $k$ level cell $C_w$ containing $x$ but not
intersecting $V_0$. ($k_0$ depends on the location of $x$ in
$\mathcal{SG}$.) Then by using Theorem 3.1 we can find a sequence of
words $w^{(k)}$ of length $k$ ($k\geq k_0$) and a sequence of mean
value neighborhoods $B_k(x)$ associated to $C_{w^{(k)}}$. Obviously,
$\{B_k(x)\}_{k\geq k_0}$ will form a system of neighborhoods of the
point $x$ satisfying $\bigcap_{k\geq k_0} B_k(x)=\{x\}$. $\Box$

\section{Mean value property of general functions on $\mathcal{SG}$}

In this section, we extend the mean value property to more general
functions on $\mathcal{SG}$. Given a point $x$ in
$\mathcal{SG}\setminus V_0$ and a cell $C_w=F_w(\mathcal{SG})$
containing $x$, for each mean value neighborhood $B$ of $x$
associated to $C_w$, we assign a constant $c_B$ to $B$. We want
$$M_B(u)-u(x)\approx c_B\Delta u(x)$$ for $u$ in $dom\Delta$.
More precisely, let $\{B_k(x)\}_{k\geq k_0}$ be the system of mean
value neighborhoods of the point $x$; we want
\begin{equation}\label{3}
\lim_{k\rightarrow\infty}\frac{1}{c_{B_k(x)}}\left(M_{B_k(x)}-u(x)\right)=\Delta
u(x)
\end{equation}
 for appropriate functions in the domain of
$\Delta$, which is the desired fractal analog of $(\ref{0})$.

 For
this purpose, let $v$ be a function on $\mathcal{SG}$ satisfying
$\Delta v\equiv 1$. For each point $x$ in $\mathcal{SG}\setminus
V_0$, and each mean value neighborhood $B$ of $x$, define $c_B$ by
$$c_B=M_B(v)-v(x).$$ Note that the result is independent of which
$v$, because any two such functions differ by a harmonic function
and the equality $M_B(h)-h(x)=0$ always holds for any harmonic
function $h$. So we can choose
$$v(x)=-\int G(x,y)d\mu(y),$$
which vanishes on the boundary of $\mathcal{SG}$. Here $G$ is
Green's function.

We will prove that $c_B$ is controlled by the size of $C_w$. More
precisely, we will prove:

 \textbf{Theorem 4.1.} \emph{Let $x\in \mathcal{SG}\setminus V_0$ and $B$ be a $k$ level mean value neighborhood of $x$. Then
$$c_0\frac{1}{5^k}\leq c_B\leq c_1\frac{1}{5^k}$$ for some constant $c_0, c_1$ which are independent of $x$.}

To prove Theorem 4.1, we need the explicit expression for the
function $v$. Recall from Section 2 that $v(x)$ is the uniform limit
of $v_M(x)$ for
$$v_M(x)=-\int G_M(x,y)d\mu(y).$$
Interchanging the integral and summation,
$$
v_M(x)=-\sum_{m=0}^M\sum_{z,z'\in V_{m+1}\setminus
V_m}g(z,z')\int\psi_{z'}^{(m+1)}(y)d\mu(y)\psi_{z}^{(m+1)}(x).
$$

Notice that for each $z\in V_{m+1}\setminus V_m$, $\psi_z^{(m+1)}$
is a piecewise harmonic spline of level $(m+1)$ satisfying
$\psi_{z}^{(m+1)}(y)=\delta_{z}(y)$ for $y\in V_{m+1}$. More
precisely, $\psi_z^{(m+1)}$ is supported in the two $(m+1)$-cells
meeting at $z$. If $F_\tau(\mathcal{SG})$ is one of these cells with
vertices $z,z_1$ and $z_2$, then
$\psi_z^{(m+1)}+\psi_{z_1}^{(m+1)}+\psi_{z_2}^{(m+1)}$ restricted to
$F_\tau(\mathcal{SG})$ is identically $1$. Thus
$$\int_{F_\tau(\mathcal{SG})}(\psi_z^{(m+1)}+\psi_{z_1}^{(m+1)}+\psi_{z_2}^{(m+1)})d\mu=\mu(F_\tau(\mathcal{SG}))=\frac{1}{3^{m+1}}.$$
By symmetry all three summands have the same integral, so
$\int_{F_\tau(\mathcal{SG})}\psi_z^{(m+1)}d\mu=\frac{1}{3^{m+2}}$.
Together with the contribution from the other $(m+1)$-cell we find
for each $z\in V_{m+1}\setminus V_m$,
\begin{equation}\label{6}
\int\psi_{z}^{(m+1)}(y)d\mu(y)=\frac{2}{3^{m+2}}.\end{equation}
Hence
$$
v_M(x)=-\frac{2}{9}\sum_{m=0}^M\frac{1}{3^m}\sum_{z,z'\in
V_{m+1}\setminus V_m}g(z,z')\psi_{z}^{(m+1)}(x).
$$
Substituting the exact value of $g(z,z')$(see Section 2 and details
in $\cite{Str2}$ page $50$) into it, we get
\begin{eqnarray*}
v_M(x)&=&-\frac{2}{9}\sum_{m=0}^M\frac{1}{3^m}\left(\sum_{|w|=m}\sum_{z,z'\in
F_w(V_0)\setminus F_w(V_1)
}g(z,z')\psi_{z}^{(m+1)}(x)\right)\\&=&-\frac{2}{9}\sum_{m=0}^M\frac{1}{3^m}\left(\sum_{|w|=m}\sum_{z\in
F_w(V_0)\setminus F_w(V_1)}(\frac{9}{50} r^m+2\frac{3}{50}r^m)\psi_{z}^{(m+1)}(x)\right)\\
&=&-\frac{1}{15}\sum_{m=0}^M\frac{1}{5^m}\phi_m(x)
\end{eqnarray*}
for $$\phi_m(x)=\sum_{z\in V_{m+1}\setminus V_m}\psi_z^{(m+1)}(x).$$
Thus
$$v(x)=-\frac{1}{15}\sum_{m=0}^{\infty}\frac{1}{5^m}\phi_m(x).$$

\textbf{Remark.} \textit{The function $v$ is invariant under
Dihedral-$3$ symmetry.}

This is a direct corollary of the fact that each $\phi_m(x)$ is
invariant under $D_3$ symmetry.

Due to the above remark, we may assume that $D_w$ associated to
$C_w$ has a fixed shape as shown in Fig. 3.1 without loss of
generality. We now show that although $c_B$ depends on the relative
position of $x$ in  $C_w$, it does not depend on the location of $x$
or $C_w$ in $\mathcal{SG}$.

\textbf{Lemma 4.1.} \textit{Let $x,x'$ be two distinct points in
$\mathcal{SG}\setminus V_0$. Let $C_w$ and $C_{w'}$ be two $k$ and
$k'$ level neighboring cells of $x$ and $x'$ respectively. Denote by
$B$ and $B'$ two mean value neighborhoods of $x$ and $x'$
respectively. If $B$ and $B'$ have the same shapes (the same
relative locations associated to $C_w$ and $C_{w'}$ respectively),
then
$$c_B=5^{k'-k}c_{B'}.$$ In particular, if $B$ and $B'$ have the same
levels and same shapes, then $c_B=c_{B'}$.}

\textit{Proof.}  $D_w$ can be decomposed into a union of a $k$ level
cell $D_{w}^{(1)}$ and a $(k-1)$ level cell $D_w^{(2)}$ as shown in
Fig. 3.1. Denote by $q$ the junction point connecting $D_{w}^{(1)}$
and $D_w^{(2)}$.  Similarly, $D_{w'}$ can also be written as a union
of a $k'$ cell  $D_{w'}^{(1)}$ and a $(k'-1)$ cell  $D_{w'}^{(2)}$
with a junction point $q'$ connecting them.

Let $\tau$ be the linear function mapping $D_w$ onto $D_{w'}$.
Suppose $D_{w}^{(1)}=F_{\alpha}(\mathcal{SG})$ and
$D_{w}^{(2)}=F_{\beta}(\mathcal{SG})$ where $\alpha$ and $\beta$ are
the corresponding words of $D_{w}^{(1)}$ and $D_{w}^{(2)}$
respectively. Similarly, denote by $\alpha'$ and $\beta'$ the
corresponding words of $D_{w'}^{(1)}$ and $D_{w'}^{(2)}$. Hence we
can write $\tau$ as $\tau(z)=F_{\alpha'}\circ F_{\alpha}^{-1}(z)$ if
$z\in D_{w}^{(1)}$, and $\tau(z)=F_{\beta'}\circ F_{\beta}^{-1}(z)$
if $z\in D_{w}^{(2)}$. In particular, $\tau(q)=q'$ and $\tau(x)=x'$.

Consider the function $(v\circ F_{\alpha}-5^{k'-k}v\circ
F_{\alpha'})$ defined on $\mathcal{SG}$. Noting that $|\alpha|=k$
and $|\alpha'|=k'$, using the scaling property of $\Delta$(see
details in $\cite{Str2}$ page $33$), we have
$$\Delta(v\circ F_{\alpha}-5^{k'-k}v\circ F_{\alpha'})=r^{|\alpha|}\frac{1}{3^{|\alpha|}}\Delta v\circ F_{\alpha}-5^{k'-k}r^{|\alpha'|}\frac{1}{3^{|\alpha'|}}\Delta v\circ F_{\alpha'}=0,$$
which shows that the difference between $v\circ F_{\alpha}$ and
$5^{k'-k}v\circ F_{\alpha'}$ is a harmonic function. Hence the
difference between $v$ and $5^{k'-k}v\circ\tau$ on $D_{w}^{(1)}$ is
 harmonic. A similar discussion will show that the difference between
$v$ and $5^{k'-k}v\circ\tau$ on $D_{w}^{(2)}$ is also harmonic.
Since the matching condition on normal derivatives of
$(v-5^{k'-k}v\circ\tau)$ at $q$ holds obviously, we have proved that
$\Delta(v-5^{k'-k}v\circ\tau)=0$ on $D_{w}$, i.e., the function
$(v-5^{k'-k}v\circ\tau)$ is harmonic on $D_{w}$.

By the definition $c_B=M_B(v)-v(x)$ and $c_{B'}=M_{B'}(v)-v(x')$.
Notice that for the second equality, by changing variables we can
write $c_{B'}=M_{B}(v\circ\tau)-v\circ\tau(x)$. Hence
$$c_B-5^{k'-k}c_{B'}=M_B(v-5^{k'-k}v\circ\tau)-(v-5^{k'-k}v\circ\tau)(x)=0,$$
since $(v-5^{k'-k}v\circ\tau)$ is a harmonic function on $D_{w}$.
$\Box$

\emph{Proof of Theorem 4.1.}

\textbf{Estimate of $c_B$ from above.}

From Lemma 4.1, since $c_B$ depends only on the relative geometry of
$B$ and $C_w$ and the size of $C_w$, but not on the location of
$C_w$, we may assume that $D_w$ is contained in a $(k-2)$ level cell
$C$ in $\mathcal{SG}$ without loss of generality.

By the definition of $c_B$, we may write
$$c_B=M_{B}(v)-v(x)=\lim_{M\rightarrow\infty}\left(\frac{1}{\mu(B)}\int_{B}v_Md\mu-v_M(x)\right).$$
Substituting the exact formula of $v_M$ into it, we get
$$c_B=-\frac{1}{15}\sum_{m=0}^\infty\frac{1}{5^m}\left(M_B(\phi_m)-\phi_m(x)\right),$$
for $$\phi_m=\sum_{z\in V_{m+1}\setminus V_m}\psi_z^{(m+1)}.$$

Notice that each $\phi_{m}$ is a piecewise harmonic spline of level
$m+1$. So when $m+1\leq k-2$, $\phi_{m}$ is harmonic in the cell
$C$, which yields that $M_{B}(\phi_m)-\phi_m(x)=0$. So the first
$k-2$ terms in the infinite series of $v$ will contribute $0$ to
$c_B$. Hence
$$c_B=-\frac{1}{15}\sum_{m=k-2}^\infty\frac{1}{5^m}\left(M_B(\phi_m)-\phi_m(x)\right).$$
It is easy to see that this implies
$$|c_B|\leq\frac{1}{15}\sum_{m=k-2}^\infty\frac{1}{5^m}\frac{1}{\mu(B)}\int_{B}|\phi_m(y)-\phi_m(x)|d\mu(y).$$
Then by the maximum principle, we finally get
$$|c_B|\leq\frac{1}{15}\sum_{m=k-2}^\infty\frac{1}{5^m}=\frac{25}{12}\cdot\frac{1}{5^k}.$$

\textbf{Estimate of $c_B$ from below.}

\begin{figure}[h]
\begin{center}
\includegraphics[width=4.6cm,totalheight=4.cm]{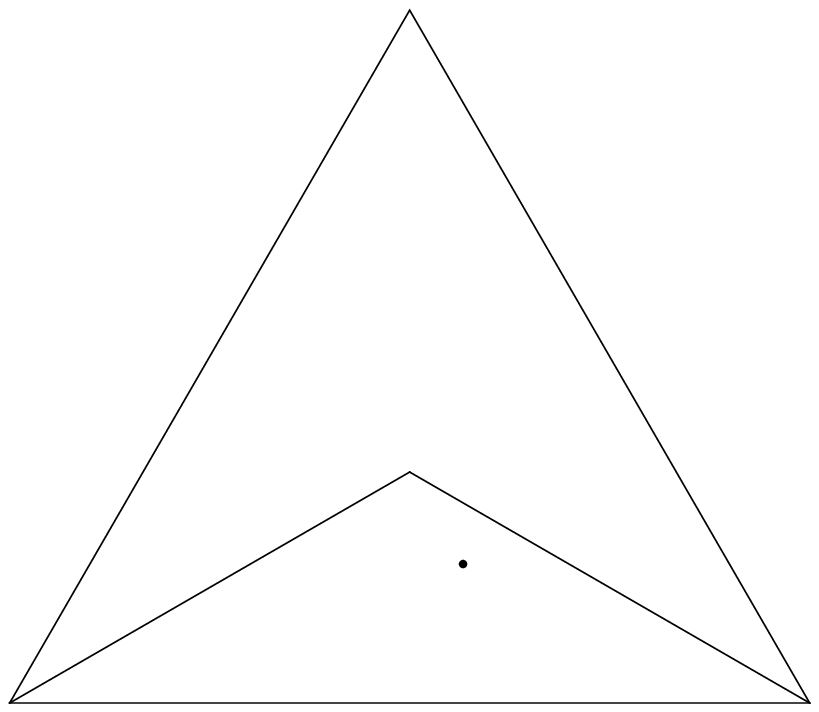}
\setlength{\unitlength}{1cm}
\begin{picture}(0,0) \thicklines
\put(-2.6,1.5){$o$}\put(-2.6,4.0){$p_0$}\put(-5.0,-0.2){$p_1$}\put(-0.4,-0.2){$p_2$}\put(-2.2,0.6){$x$}\put(-1.8,1.3){$C_w$}
\end{picture}
\begin{center}
\textbf{Fig. 4.1.} \small{a $1/3$ region of $C_w$.}
\end{center}
\end{center}
\end{figure}
Without loss of generality, we assume that $x$ is located in the
$1/3$ region of $C_w$ as shown in Fig. 4.1, i.e., $x$ is contained
in the triangle $T_{p_1,p_2,o}$, where $o$ is the geometric center
of $C_w$. Then by the proof of Theorem 3.1, $B$ is a subset of the
union of $C_w$ and two of its neighbors $C_1$ and $C_2$. Hence we
can write $B=C_w\cup E_1\cup E_2$, where $E_i=B\cap C_i$.

\textbf{Claim 1.}\emph{ Let $\widetilde{B}=F_0(\mathcal{SG})\cup
\widetilde{E}_1\cup \widetilde{E}_2$, where $\widetilde{E}_i$ is a
triangle obtained by cutting $F_i(\mathcal{SG})$ symmetrically with
a line below the top vertex $F_iq_0$.(see Fig. 4.2.) If
$\widetilde{B}$ and $B$ have the same shapes, then
$$c_B=5^{1-k}c_{\widetilde{B}}.$$}\begin{figure}[h]
\begin{center}
\includegraphics[width=4.6cm,totalheight=4cm]{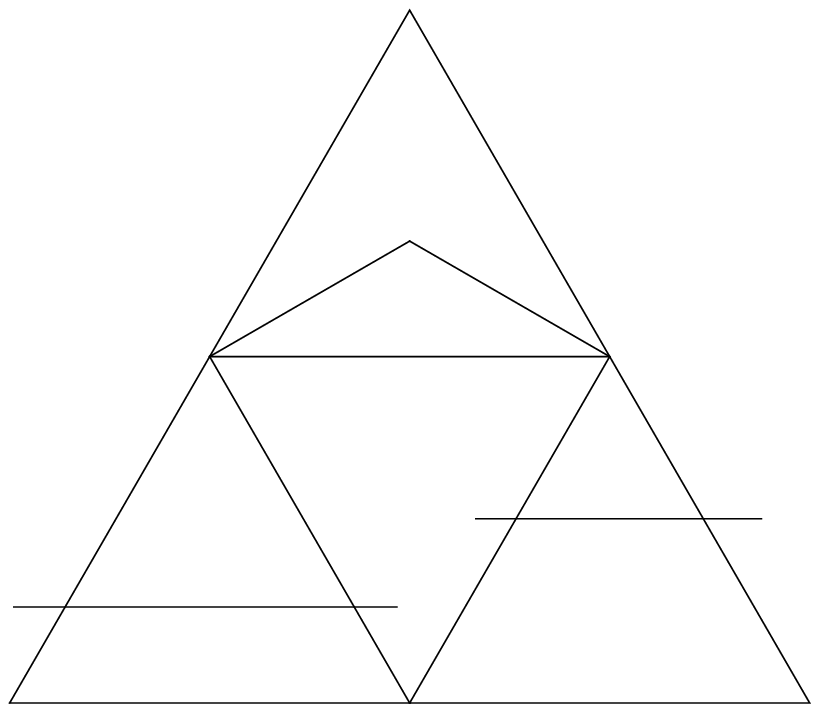}
\setlength{\unitlength}{1cm}
\begin{picture}(0,0) \thicklines
\put(-3.2,2.8){$F_0(\mathcal{SG})$}\put(-2.6,4.0){$q_0$}\put(-5.0,-0.2){$q_1$}\put(-0.4,-0.2){$q_2$}
\put(-1.7,1.2){$\widetilde{E}_2$}\put(-3.9,0.8){$\widetilde{E}_1$}
\end{picture}
\begin{center}
\textbf{Fig. 4.2.} \small{a sketch of $\widetilde{B}$.}
\end{center}
\end{center}
\end{figure}

This is a direct corollary of Lemma 4.1.

We only need to prove that $c_{\widetilde{B}}$ for $\widetilde{B}$
defined in Claim 1 has a positive lower bound. For simplicity of
notation, in all that follows, we write $B$ instead of
$\widetilde{B}$. In other words, we only need to consider $B$ whose
associate cell $C_w$ is $F_{0}(\mathcal{SG})$. In this setting,
$p_i=F_0q_i$, $C_1=F_1(\mathcal{SG})$ and $C_2=F_2(\mathcal{SG})$.

We write $v=-\frac{1}{15}\widetilde{v}$ where $\widetilde{v}$ is the
non-negative function defined by
$$\widetilde{v}=\sum_{m=0}^{\infty}\frac{1}{5^m}\phi_m.$$  For each $M\geq
0$, denote by
$$\widetilde{v}_M=\sum_{m=0}^{M}\frac{1}{5^m}\phi_m$$ the partial sum of the first $M+1$ terms of
$\widetilde{v}$. Then $\widetilde{v}_M$ converges to $\widetilde{v}$
uniformly as $M\rightarrow \infty$.

We have the following three claims on $\widetilde{v}$.

\textbf{Claim 2.} \emph{$0\leq \widetilde{v}\leq 1$ on
$\mathcal{SG}$ and $\widetilde{v}$ takes constant $1$ along the
maximal inner upside-down triangle contained in $\mathcal{SG}$.}

\emph{Proof.}  Consider the partial sum function $\widetilde{v}_M$.
Obviously, $\widetilde{v}_M$ is a $(M+1)$-level piecewise harmonic
function on $\mathcal{SG}$. For convenience, denote by $\nabla$ the
maximal inner upside-down triangle contained in $\mathcal{SG}$. We
divide the vertices $V_{M+1}$ into three parts, $V'_{M+1}$,
$V''_{M+1}$ and $V'''_{M+1}$, where $V'_{M+1}$ consists of those
vertices lying along $\nabla$, $V''_{M+1}$ consists of those
vertices at distance $2^{-(M+1)}$ from $\nabla$, and $V'''_{M+1}$
consists of the remain vertices. Then by using the
$``\frac{1}{5}-\frac{2}{5}"$ rule, an inductive argument shows that
$\widetilde{v}_M\equiv1$ on $V'_{M+1}$,
$\widetilde{v}_M\equiv1-\frac{1}{5^{M}}$ on $V''_{M+1}$, and
$\widetilde{v}_M\leq 1-\frac{1}{5^M}$ on $V'''_{M+1}$. Since
$\widetilde{v}$ is the uniform limit of $\widetilde{v}_M$ and
$V'_{M+1}$ goes to $\nabla$ as $M$ goes to the infinity, we then
have $0\leq \widetilde{v}\leq 1$ on $\mathcal{SG}$ and
$\widetilde{v}\equiv 1$ on $\nabla$. $\Box$

\textbf{Claim 3.} \emph{For each $x$ contained in the triangle
$T_{p_1,p_2,o},$ $\widetilde{v}(x)\geq\frac{24}{25}.$}

\emph{Proof.}  For $\tau=(0,1,1), (0,1,2), (0,2,1)$ and $(0,2,2)$,
by using the harmonic extension algorithm, namely, the
$``\frac{1}{5}-\frac{2}{5}"$ rule, we get that
$$
\widetilde{v}(F_{\tau}q_0)=\widetilde{v}_2(F_{\tau}q_0)=\sum_{m=0}^{2}\frac{1}{5^m}\phi_m(F_\tau
q_0)=1\cdot\frac{4}{5}+\frac{1}{5}\cdot\frac{3}{5}+\frac{1}{25}\cdot
1=\frac{24}{25},$$ where $\frac{4}{5}$, $\frac{3}{5}$ and $1$ are
the values of $\phi_0$, $\phi_1$ and $\phi_2$ at $F_\tau q_0$
respectively. Also, for those $\tau$, by Claim 2, we have
$$\widetilde{v}(F_{\tau}q_1)=\widetilde{v}_2(F_{\tau}q_1)=\widetilde{v}(F_{\tau}q_2)=\widetilde{v}_2(F_{\tau}q_2)=1.$$
  Notice that for each point $x$ in the
triangle $T_{p_1,p_2,o}$, $x$ is contained in one of the four $3$-
level cells $F_{011}(\mathcal{SG})$, $F_{012}(\mathcal{SG})$,
$F_{021}(\mathcal{SG})$ and $F_{022}(\mathcal{SG})$. Since
$\widetilde{v}_2$ is harmonic in each such cell, by using the
maximal principle, we get that
$$\widetilde{v}_2(x)\geq \frac{24}{25}.$$ Hence
$\widetilde{v}(x)\geq \frac{24}{25}$ since each term in the infinite
series of $\widetilde{v}$ is non-negative. $\Box$

\textbf{Claim 4.} \emph{$M_B(\widetilde{v})\leq\frac{17}{18}.$}

\emph{Proof.} First of all we prove that
$$\int_{F_0(\mathcal{SG})}\widetilde{v}(y)d\mu(y)=\frac{5}{18}.$$ We need to compute
$\int_{F_0(\mathcal{SG})}\phi_{m}(y)d\mu(y)$ for each non-negative
integer $m$. For each $m\geq 0$,
$$\int_{F_0(\mathcal{SG})}\phi_{m}(y)d\mu(y)=\frac{1}{3}\cdot 3^{m+1}\cdot\frac{2}{3^{m+2}}=\frac{2}{9},$$
by using $(\ref{6})$ and the fact that $\phi_m=\sum_{z\in
V_{m+1}\setminus V_m} \psi_{z}^{(m+1)}$.
 Hence
$$\int_{F_0(\mathcal{SG})}\widetilde{v}(y)d\mu(y)=\frac{2}{9}\sum_{m=0}^{\infty}\frac{1}{5^m}=\frac{5}{18}.$$

By our assumption, the mean value neighborhood $B$ can be written as
$$B=F_0(\mathcal{SG})\cup E_1\cup E_2,$$ where $E_i=B\cap C_i$.   Hence we have
\begin{eqnarray*}
M_B(\widetilde{v})
&=&\frac{1}{\mu(B)}\left(\int_{F_0(\mathcal{SG})}\widetilde{v}(y)d\mu(y)+\int_{E_1}\widetilde{v}(y)d\mu(y)+\int_{E_2}\widetilde{v}(y)d\mu(y)\right)\\
&\leq
&\frac{1}{\mu(B)}\left(\int_{F_0(\mathcal{SG})}\widetilde{v}(y)d\mu(y)+\int_{E_1}1\cdot
d\mu(y)+\int_{E_2}1\cdot d\mu(y)\right)\\
&=&\frac{{5}/{18}+\mu(E_1)+\mu(E_2)}{{1}/{3}+\mu(E_1)+\mu(E_2)},
\end{eqnarray*}
where the inequality follows from Claim 2.  Since $0\leq
\mu(E_1)+\mu(E_2)\leq \frac{2}{3}$,
$\frac{\frac{5}{18}+x}{\frac{1}{3}+x}$ is increasing in $x\geq 0$,
$$\frac{{5}/{18}+\mu(E_1)+\mu(E_2)}{{1}/{3}+\mu(E_1)+\mu(E_2)}\leq\frac{5/18+2/3}{1/3+2/3}=\frac{17}{18}.$$
Hence we always have
$$M_B(\widetilde{v})\leq\frac{17}{18}.\quad \Box$$

Now we turn to estimate $c_B$. Obviously,
\begin{equation*}
c_B=M_B(v)-v(x)=-\frac{1}{15}\left(M_B(\widetilde{v})-\widetilde{v}(x)\right).
\end{equation*}

By Claim 3 and Claim 4, we notice that
$M_B(\widetilde{v})-\widetilde{v}(x)\leq\frac{17}{18}-\frac{24}{25}=-\frac{7}{450}$.
Hence
$$c_B\geq\frac{1}{15}\cdot\frac{7}{450}>0.\quad \Box$$

On the other hand, given a point $x$ and $C_w=F_w(\mathcal{SG})$ a
$k$ level neighborhood of $x$, for any $u\in dom\Delta$, we write
$$u=h^{(k)}+(\Delta u(x))v+R^{(k)}$$ on $C_w$,
where $h^{(k)}$ is a harmonic function defined by
$$h^{(k)}+(\Delta u(x))v|_{\partial C_w}=u|_{\partial C_w}.$$

It is not hard to prove the following estimate:

\textbf{Lemma 4.2.}\emph{ Let $u\in dom\Delta$ with $g=\Delta u$
satisfying the following H\"{o}lder condition
$$|g(y)-g(x)|\leq c\gamma ^k, \quad \quad (0<\gamma<1)$$ for all $y\in C_w$. Then the remainder satisfies
$$R^{(k)}=O\left((\frac{\gamma}{5})^k\right)$$ on $C_w$(hence also on
$B_k(x)$).}

\emph{Proof.} It is easy to check that $\Delta R^{(k)}(y)=\Delta
u(y)-\Delta u(x)$ and $R^{(k)}(y)$ vanishes on the boundary of
$C_w$. Hence $R^{(k)}$ is given by the integral of $\Delta
u(y)-\Delta u(x)$ on $C_w$ against a scaled Green's function.
Noticing that the scaling factor is $(\frac{1}{5})^k$ and
$$|\Delta u(y)-\Delta u(x)|\leq c\gamma^k,$$ we then get the desired result. $\Box$

This looks like a Taylor expansion remainder estimate of $u$ at $x$.
See more details on this topic in $\cite{Str3}$.

\textbf{Remark.} \emph{If we require $u\in dom \Delta^2$, then the
remainder $R^{(k)}$ satisfies
$$R^{(k)}=O\left((\frac{3}{5}\cdot\frac{1}{5})^k\right)$$ on $C_w$(hence also on
$B_k(x)$). The reason is that in this case $\Delta u$ satisfies the
H\"{o}lder condition that $|\Delta u(y)-\Delta u(x)|\leq
c(\frac{3}{5})^k$ for all $y\in C_w$, because $\Delta^2 u$ is
assumed continuous, see $\cite{Str3}$, Theorem 8.4.}

Using the above lemma and Theorem 4.1, we then have the following
main result of this section.

\textbf{Theorem 4.2.} \emph{Let $u\in dom\Delta$ with $g=\Delta u$
satisfying the H\"{o}lder condition $|g(y)-g(x)|\leq c\gamma ^k$ for
some $\gamma$ with $ 0<\gamma<1$, for all $x,y$ belonging to the
same $k$ level cell. Then
$$\lim_{k\rightarrow\infty}\frac{1}{c_{B_k(x)}}\left(M_{B_{k}(x)}(u)-u(x)\right)=\Delta u(x).$$}

\emph{Proof.}  Using Taylor expansion of $u$ and noticing that
$M_{B_{k}(x)}(h^{(k)})-h^{(k)}(x)=0$,
$M_{B_{k}(x)}(v)-v(x)=c_{B_k(x)}$, we have
\begin{eqnarray*}
\frac{1}{c_{B_k(x)}}\left(M_{B_{k}(x)}(u)-u(x)\right)-\Delta
u(x)&=&\frac{1}{c_{B_k(x)}}\left(M_{B_{k}(x)}(R^{(k)})-R^{(k)}(x)\right)\\
&=&\frac{1}{c_{B_k(x)}}O\left((\frac{\gamma}{5})^k\right)=O(\gamma^k).
\end{eqnarray*}
Hence letting $k\rightarrow\infty$, we get the desired result.
$\Box$

\section{p.c.f. fractals with Dihedral-3 symmetry}

The results for $\mathcal{SG}$ should extend to other p.c.f.
fractals which possess symmetric properties of both the geometric
structure and the harmonic structure. We assume that a regular
harmonic structure is given on a p.c.f. self-similar fractal $K$.
The reader is referred to $\cite{Ki2,Str2}$ for exact definitions
and any unexplained notations.  We assume now that $\sharp V_0=3$
and all structures possess full $D_3$ symmetry. This means there
exists a group $\mathcal{G}$ of homeomorphisms of $K$ isomorphic to
$D_3$ that acts as permutations on $V_0$, and $\mathcal{G}$
preserves the harmonic structures and the self-similar measure.

Assume that the fractal $K$ is the invariant set of a finite
iterated function system of contractive similarities. We denote
these maps $\{F_i\}_{i=1,\cdots,N}$ with $N\geq 3$. Let $r_i$ denote
the $i$-th resistance renormalization factor and $\mu_i$ denote the
$i$-th weight of the self-similar measure $\mu$ on $K$. In general,
it is not necessary that all $r_i$'s and all $\mu_i$'s be the same,
but here we must have $r_0=r_1=\cdots=r_N$ and
$\mu_0=\mu_1=\cdots=\mu_N$ from the above Dihedral-3 symmetry
assumption. We denote $V_0=\{q_0,q_1,q_2\}$ the set of boundary
points.

\textbf{Examples.} (i) The Sierpinski gasket $\mathcal{SG}$. In this
case all $r_i=3/5$ and all $\mu_i=1/3$.

(ii) The hexagasket, or fractal Star of David, can be generated by
$6$ maps with simultaneously rotate and contract by a factor of
$1/3$ in the plane. Thus $V_0$ consists of $3$ points of an
equilateral triangle, and $V_1$ consists of the vertices of the Star
of David, as shown in Fig. 5.1. Although the same geometric fractal
can be constructed by using contractions which do not rotate, this
gives rise to a different self-similar structure (in particular with
$\sharp V_0=6$). Our choice of self-similar structure destroys the
$D_6$ symmetry of the geometric fractal, but it has the advantage of
easier computation. In this case, all $r_i=3/7$ and all $\mu_i=1/6$.
Note that in this example there exist points in $V_1$ that are not
junction points.
\begin{figure}[h]
\begin{center}
\includegraphics[width=10cm,totalheight=5.cm]{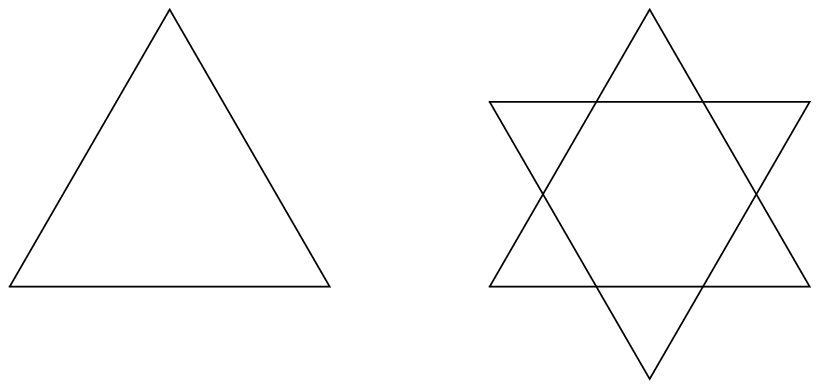}
\setlength{\unitlength}{1cm}
\begin{picture}(0,0) \thicklines
\put(-8.4,-0.2){$\Gamma_0$}\put(-2.5,-0.2){$\Gamma_1$}
\end{picture}
\begin{center}
\textbf{Fig. 5.1.} \small{The first $2$ graphs, $\Gamma_0,\Gamma_1$
in the approximation to the hexagasket.}
\end{center}
\end{center}
\end{figure}
(iii) The level $3$ Sierpinski gasket $\mathcal{SG}_3$, obtained by
taking $6$ contractions of ratio $1/3$ as shown in Fig. 5.2. Here we
have all $r_i=7/15$ and $\mu_i=1/6$. Note that all seven vertices in
$V_1\setminus V_0$ are junction points, but the one in the middle
intersects three $1$-cells. In a similar manner we could define
$\mathcal{SG}_n$ for any value of $n\geq 2$.
\begin{figure}[h]
\begin{center}
\includegraphics[width=6cm,totalheight=5.cm]{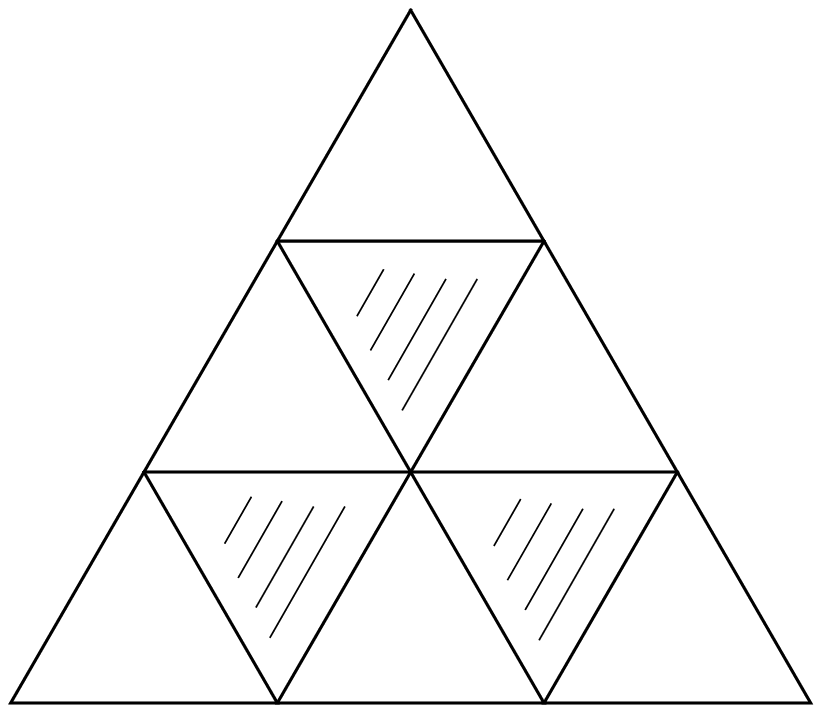}
\begin{center}
\textbf{Fig. 5.2.} \small{The graph of the $V_1$ vertices of the
level $3$ Sierpinski gasket.}
\end{center}
\end{center}
\end{figure}

We prove that there are results analogous to Theorem 3.1, which
yield the existence of mean value neighborhoods associated to $K$.

Given a point $x$ in $K\setminus V_0$, consider any cell $F_wK=C_w$
with boundary points $p_0,p_1,p_2$ containing the point $x$. Without
losing of generality, we may require that the cell $C_w$ does not
intersect $V_0$. For each $i$, denote by $C_{i,1},\cdots,C_{i,l_i}$
the neighboring cells of $C_w$ of the same size, intersecting $C_w$
at $p_i$, where $l_i$ is the number of such cells. It is possible
that $l_i=0$ for some $i$ since $p_i$ may be a non-junction point.
If this is true, the matching condition says that the normal
derivative of any harmonic function $h$ must be zero at this point,
which yields that the value of $h$ at this point is the mean value
of the values of $h$ at the other two boundary points of $C_w$. In
other words, the restriction of all global harmonic functions in
$C_w$  is two dimensional. Denote by $D_w$ the union of $C_w$ and
all its neighboring cells, i.e.,
$$D_w=C_w\cup\bigcup_{i,j}C_{i,j}.$$ Two cells $C_w$ and $C_{w'}$
are said to have the same \emph{neighborhood type} if they have the
same relative geometry with respect to $D_w$ and $D_{w'}$
respectively. It is obvious that there only exist finitely many
distinct types. For example, for $\mathcal{SG}$, all cells have
exactly only one neighborhood type. For $\mathcal{SG}_3$, the number
of the finite types is $3$. For $\mathcal{SG}_n (n\geq 4)$, the
number of the finite types becomes $4$. For the hexagasket gasket,
the number of the finite types is $2$.

Let $h$ be a harmonic function on $K$. Given a set $B$ containing
$C_w$, define
$$M_B(h)=\frac{1}{\mu(B)}\int_B h d\mu$$ the \emph{mean value} of $h$ over $B$. We
are interested in an identity
\begin{equation}\label{4}
M_B(h)=\sum_i a_i h(p_i)
\end{equation}
for some coefficients $(a_0,a_1,a_2)$ satisfying $\sum a_i=1$.
Notice that this is true for $\mathcal{SG}$. In that setting, a
harmonic function is uniquely determined by its values on the
boundary of any given cell $C_w$ because the harmonic extension
matrix associated with $C_w$ is invertible. However, in the general
case, the harmonic extension matrices may not be invertible. So we
can not prove $(\ref{4})$ for every set $B$ simply by linearity.
However, it will suffice to show that the equality $(\ref{4})$ holds
for certain specified sets $B$.

Consider a set $B$ which is a subset of $D_w$, containing $C_w$.
Then $B$ must be made up of four parts, i.e.,
$$B=C_w\cup E_0\cup E_1\cup E_2$$ where $E_i=B\cap C_i$ with $C_i=\cup_{j=1}^{l_i}C_{i,j}$.
 It is possible that $C_i$ may be  empty  since $p_i$
may be a nonjunction point. We can also subdivide each $E_i$ into
$l_i$ small pieces, i.e., $E_i=\cup_j E_{i,j}$ for $E_{i,j}=E_i\cap
C_{i,j}$. For each $i$, we require that $E_{i,1},\cdots,E_{i,l_i}$
be of the same size and shape. Moreover, in analogy with the
$\mathcal{SG}$ case, we require that each $E_{i,j}$ to be a
symmetric (under the reflection symmetry that fixes $p_i$) cutoff
sub-triangle of $C_{i,j}$, containing $p_i$ as one of its vertex
points. This means that there is a straight line $L_{i,j}$,
symmetric under the reflection symmetry fixing $p_i$, cutting
$C_{i,j}$ into two parts, and $E_{i,j}$ is the one containing $p_i$.
For each $E_{i,j}$, define the distance between $p_i$ and the line
$L_{i,j}$ the size of $E_{i,j}$. Of course, for each fixed $i$,
$E_{i,1},\cdots,E_{i,l_i}$ have the same sizes. We call the common
value the size of $E_i$. Suppose the size of every $C_{i,j}$ is
$\rho$. (Of course, they are all equal.) Then for each $i$, the size
of $E_i$ is $c_i \rho$ where the coefficient $0 \leq c_i\leq 1$.
Hence we can write the set $B=B(c_0,c_1,c_2)$. (If $p_i$ is a
nonjunction point, then $c_i$ should always be $0$.) For example,
suppose that the boundary points of $C_w$ consist of junction
points, then $B(0,0,0)=C_w$ and $B(1,1,1)=D_w$. Denote by
$$\mathcal{B}=\{B(c_0,c_1,c_2): 0\leq c_i\leq 1\}$$ the family of
all such sets. Then we can show that the formula $(\ref{4})$ holds
for each $B\in \mathcal{B}$.

\textbf{Proposition 5.1.} \emph{Let $B\in \mathcal{B}$, then for any
harmonic function $h$, we have $(\ref{4})$ for some coefficients
$(a_0,a_1,a_2)$ independent of $h$. Moreover, $\sum_i a_i=1$.}

\textit{Proof.} Each $B\in \mathcal{B}$ can be written as $B=C_w\cup
E_0\cup E_1 \cup E_2$. Given a harmonic function $h$ on $K$, for
fixed $i$, we first consider the integral $\int_{E_i}hd\mu$.
Obviously,
$$\int_{E_i}hd\mu=\sum_j \int_{E_{i,j}}hd\mu.$$ For each $1\leq j\leq l_i$, denote by
$\{z_{i,j},w_{i,j},p_i\}$ the boundary points of $C_{i,j}$. Since
each $E_{i,j}$ is contained in $C_{i,j}$,
$\frac{1}{\mu(C_w)}\int_{E_{i,j}}hd\mu$ can be expressed as a linear
combination of $h(p_i),h(z_{i,j})$ and $h(w_{i,j})$ with
non-negative  coefficients independent of the harmonic function $h$.
Since the set $E_{i,j}$ is symmetric under the reflection symmetry
fixing $p_i$, the two  coefficients with respect to $h(z_{i,j})$ and
$h(w_{i,j})$ must be equal. In other words, we can write
$$\int_{E_{i,j}}hd\mu=\left(m_{i,j}h(p_i)+n_{i,j}h(z_{i,j})+n_{i,j}h(w_{i,j})\right)\mu(C_w)$$
for $m_{i,j},n_{i,j}\geq 0$. Moreover, since for each fixed $i$,
$E_{i,j}$ are in the same relative position associated to $C_{i,j}$
for different $j$'s, $\int_{E_{i,j}}hd\mu$ can be expressed as a
linear combination of $h(p_i),h(z_{i,j}), h(w_{i,j})$ with the same
coefficients for different $j$'s. Hence we can write
$$\int_{E_i}hd\mu=\left(m_ih(p_i)+n_{i}\sum_{j}\left(h(z_{i,j})+h(w_{i,j})\right)\right)\mu(C_w),$$
for suitable coefficients $m_i,n_i\geq 0$. The mean value property
at the point $p_i$ says that
$$\sum_{j}(h(z_{i,j})+h(w_{i,j}))=(2l_i+2)h(p_i)-(h(p_{i-1})+h(p_{i+1})).$$
Combining the above two equalities, we get
$$\int_{E_i}hd\mu=\left((m_i+2l_in_i+2n_i)h(p_i)-n_{i}h(p_{i-1})-n_{i}h(p_{i+1})\right)\mu(C_w).$$

On the other hand, by the linearities and symmetries of both the
harmonic structure and the self-similar measure,
$$\int_{C_w}hd\mu=\frac{\mu(C_w)}{3}\left(h(p_0)+h(p_1)+h(p_2)\right).$$

Since the ratio of $\mu(E_{i,j})$ to $\mu(C_w)$ depends only on
$c_i$, we have proved that $M_B(h)$ can be viewed as a linear
combination of the values of $h$ on the boundary points of $C_w$,
i.e.,
$$M_B(h)=\sum_i a_i h(p_i),$$ where the combination coefficients are
independent of $h$. Moreover, we must have $\sum a_i=1$ by
considering $h\equiv 1.$
 $\Box$

\textbf{Remark 1.} \textit{This means that $M_B(h)$ is a weighted
average of the values $h(p_0), h(p_1)$ and $h(p_2)$. Moreover, if
one of the boundary points, for example $p_2$, is a nonjunction
point, then by the fact that $h(p_2)=\frac{1}{2}(h(p_0)+h(p_1))$, we
have
$$M_B(h)=a_0h(p_0)+a_1h(p_1)+\frac{1}{2}a_2\left(h(p_0)+h(p_1)\right)=\widetilde{a}_0h(p_0)+\widetilde{a}_1h(p_1)$$
for $\widetilde{a}_0=a_0+\frac{1}{2}a_2$ and
$\widetilde{a}_1=a_1+\frac{1}{2}a_2$. We also have
$\widetilde{a}_0+\widetilde{a}_1=1$. Hence in this case, we can also
view $M_B(h)$ as a weighted average of the values of $h(p_0)$ and
$h(p_1)$. }

\textbf{Remark 2.} \emph{The proof of Proposition 5.1 shows that
$(a_0,a_1,a_2)$ depends only on the neighborhood type of $C_w$ and
the relative position of $B$ associated to $C_w$, and does not
depend on the particular choice of $C_w$. In other words, if we
consider a cell $C_w$ with a given neighborhood type, then for each
set $B\in \mathcal{B}$ with the expression $B=B(c_0,c_1,c_2)$, the
coefficients $(a_0,a_1,a_2)$ of $B$ depend only on $(c_0,c_1,c_2)$.}

The following is the main result in this section.

\textbf{ Theorem 5.1.} \emph{Given a point $x\in K\setminus V_0$,
let $C_w$ be a cell containing $x$, not intersecting $V_0$, and let
$D_w$ be the union of $C_w$ and its neighboring cells of the same
size. Then there exists a mean value neighborhood $B$ of $x$
satisfying $C_w\subset B\subset D_w$.  Moreover, for each point
$x\in K\setminus V_0$, there exists a system of mean value
neighborhoods $B_k(x)$ with $\bigcap_k B_k(x)=\{x\}$.}

\emph{Proof.} We need to classify the distinct neighborhood types
into three cases according to the number of nonjunction points in
the set of boundary points of $C_w$.

\textbf{Case 1.} \textbf{All boundary points of $C_w$ are junction
points.}

This case is similar to what we have described in the $\mathcal{SG}$
setting. Let $W$ denote the triangle in $\mathbb{R}^3$ with boundary
points $P_0=(1,0,0), P_1=(0,1,0)$ and $P_2=(0,0,1)$ and $\pi_W$ the
plane containing $W$. Notice that from Proposition 5.1,
$(a_0,a_1,a_2)\in \pi_W$ for each $B$. We use $\mathcal{T}$ to
denote the map from $\mathcal{B}$ to $\pi_W$. From Remark 2 of
Proposition 5.1, the map $\mathcal{T}$ is uniquely determined by the
neighborhood type of $C_w$.  Let $\mathcal{B}^*$ be a subfamily
 contained in $\mathcal{B}$ defined by
$$\mathcal{B}^*=\{B(0,c_1,c_2):0\leq c_i\leq 1\}\cup\{B(c_0,0,c_2):0\leq c_i\leq 1\}\cup \{B(c_0,c_1,0): 0\leq c_i\leq
1\},$$ i.e.,  those elements $B$ in $\mathcal{B}$ which have the
decomposition form $B=C_w\cup E_1\cup E_2$ or $B=C_w\cup E_0\cup
E_2$, or $B=C_w\cup E_0\cup E_1$. Then we have

\textbf{Claim 1.} \emph{The map $\mathcal{T}$ from $\mathcal{B}$ to
$\pi_W$ fills out a region $\widetilde{W}$ which contains the
triangle $W$. Moreover, $\mathcal{T}$ is one-to-one from
$\mathcal{B}^*$ onto $\widetilde{W}$}.

\emph{Proof.} The proof is similar to the $\mathcal{SG}$ case. The
only difference is the line segments $\overline{OQ_0}$ and
$\overline{OP_2}$ described in the proof of Theorem 3.1 may become
continuous curves $\widehat{OQ_0}$ and $\widehat{OP_2}$ in the
general setting. $\Box$

\textbf{Case 2.} \textbf{There is one nonjunction point (for
example, $p_2$) among the boundary points of $C_w$.}

In this case, there is no neighboring cell intersecting $C_w$ at the
point $p_2$. Hence $E_2$ will always be empty. So
$\mathcal{B}=\{B(c_0,c_1,0):0\leq c_i\leq 1\}$ for this case.

 As shown in Remark 1 of Proposition 5.1, for any harmonic function $h$ on $K$, $B\in
 \mathcal{B}$,
 $M_B(h)$ is a weighted average of $h(p_0)$ and $h(p_1)$, i.e.,
 $$M_B(h)=a_0h(p_0)+a_1h(p_1)$$ with $a_0, a_1$ independent of
 $h$, satisfying $a_0+a_1=1$. Let
$I$ denote the line segment in $\mathbb{R}^2$ with endpoints
$P_0=(1,0), P_1=(0,1)$ and $\rho_I$ the line containing $I$. Notice
that from Remark 1 of Proposition 5.1, $(a_0,a_1)\in \rho_I$ for
each $B$. We still use $\mathcal{T}$ to denote the map from
$\mathcal{B}$ to $\rho_I$. From Remark 2 of Proposition 5.1, the map
$\mathcal{T}$ is uniquely determined by the neighborhood type of
$C_w$. We may write $\mathcal{T}(B(c_0,c_1,0))=(a_0,a_1)$
 for each set $B(c_0,c_1,0)$. We
will show the image of the map $\mathcal{T}$  covers the line
segment $I$. Similar to Case 1, let $\mathcal{B}^*$ be a subfamily
 contained in $\mathcal{B}$ defined by
$$\mathcal{B}^*=\{B(c_0,0,0):0\leq c_0\leq 1\}\cup\{B(0,c_1,0):0\leq c_1\leq 1\},$$ i.e.,  those elements $B$ in $\mathcal{B}$ which have
the decomposition form $B=C_w\cup E_0$ or $B=C_w\cup E_1$. Then we
have

\textbf{Claim 2.} \emph{The map $\mathcal{T}$ from $\mathcal{B}$ to
$\rho_I$ fills out the line segment $I$. Moreover, $\mathcal{T}$ is
a one-to-one map on $\mathcal{B}^*$.}

\emph{Proof.} The proof is similar to Case 1. Denote by
$O=(\frac{1}{2},\frac{1}{2})$ the midpoint of $I$.  We only prove
the map $\mathcal{T}$ from $\mathcal{B}$ to $\rho_I$ fills out half
of the line segment $I$. Then we will get the desired result by
symmetry.

Let $h$ be a harmonic function on $K$. We consider
$\mathcal{T}(\{(B(c,0,0)):0\leq c\leq 1\})$. When $c=0$,
$B(0,0,0)=C_w$ and $M_{C_w}(h)=\frac{1}{3}(h(p_0)+h(p_1)+h(p_2))$.
Combining this with the fact that
$$h(p_2)=\frac{1}{2}(h(p_0)+h(p_1)),$$ we get
$$M_{C_w}(h)=\frac{1}{2}(h(p_0)+h(p_1)).$$ Hence $\mathcal{T}(B(0,0,0))$ is the
midpoint $O$ of $I$. When $c=1$, $B(1,0,0)=C_w\cup C_0$, and an easy
calculation  gives that $M_{C_w\cup C_0}=h(p_0)$. Hence
$\mathcal{T}(B(1,0,0))$ is the endpoint $P_0$. So if we vary $c$
continuously  between $0$ and $1$, we can fill out the line segment
joining $O$ and $P_0$, which is half of $I$. $\Box$

 \textbf{Case 3.} \textbf{There are two nonjunction points (for example, $p_1$ and $p_2$)
among the boundary points of $C_w$.}

In this case, let $h$ be any harmonic function on $K$. By the
matching condition on both points $p_1$ and $p_2$, $h$ must be
constant on the whole cell $C_w$. Hence for every point $x\in C_w$,
we could view $C_w$ itself as the mean value neighborhood of $x$.

Hence the proof of Theorem 5.1 is completed by using a same argument
as that of Theorem 3.2. $\Box$

We should mention here that the result can also be extended to some
other p.c.f. fractals including the $3$-dimensional Sierpinski
gasket. However, it seems that some strong symmetric conditions of
both the geometric and the harmonic structures should be required.

\textbf{Acknowledgements.} This work was done while the first author
was visiting the Department of Mathematics, Cornell University. He
expresses his sincere gratitude to the department for their
hospitality. We would also like to thank the anonymous referees for
several important suggestions which led to the improvement of the
manuscript.

(Hua Qiu) DEPARTMENT OF MATHEMATICS, NANJING UNIVERISITY, NANJING,
210093, CHINA

\emph{E-mail address}: huatony@gmail.com

(Robert S. Strichartz) DEPARTMENT OF MATHEMATICS, MALOTT HALL,
CORNELL UNIVERSITY, ITHACA, NY 14853

\emph{E-mail address}: str@math.cornell.edu

\end{document}